\documentclass[11pt]{article}
\usepackage{latexsym,amssymb}
\def\demo{\noindent{\bf Proof. }}
\def\QED{\hfill$\Box$}
\newtheorem{Theorem}{Theorem}[section]
\newtheorem{Lemma}[Theorem]{Lemma}
\newtheorem{Corollary}[Theorem]{Corollary}
\newtheorem{Proposition}[Theorem]{Proposition}
\newtheorem{Example}[Theorem]{Example}
\newtheorem{Conjecture}[Theorem]{Conjecture}
\newtheorem{Definition}[Theorem]{Definition}

\begin{document}
\topmargin3mm
\hoffset=-1cm
\voffset=-1.5cm
\

\medskip

\begin{center}
{\large\bf Cohen-Macaulay clutters with 
combinatorial optimization properties and parallelizations of 
normal edge ideals 
}
\vspace{6mm}\\

\medskip

Luis A. Dupont, Rafael H. Villarreal\footnote{Partially supported by
CONACyT  
grant 49251-F and SNI.} 
\\ 
{\small Departamento de Matem\'aticas}\vspace{-1mm}\\ 
{\small Centro de Investigaci\'on y de Estudios Avanzados del
IPN}\vspace{-1mm}\\   
{\small Apartado Postal 14--740}\vspace{-1mm}\\ 
{\small 07000 M\'exico City, D.F.}\vspace{-1mm}\\ 
{\small e-mail: {\tt vila@math.cinvestav.mx}}\vspace{4mm}

Enrique Reyes\footnote{Partially supported by
COFAA-IPN and SNI.} 
\\
{\small Departamento de Ciencias B\'asicas, Unidad
Profesional}\vspace{-1mm}\\ 
{\small Interdisciplinaria en Ingenier\'\i a y Tecnologias Avanzadas
del IPN, }\vspace{-1mm}\\  
{\small UPIITA, Av. IPN 2580, Col. Barrio la Laguna
Ticom\'an}\vspace{-1mm}\\   
{\small 07340 M\'exico City, D.F.}\vspace{-1mm}\\

\end{center}
\date{}

\begin{abstract} 
\noindent Let $\mathcal{C}$ be a uniform clutter and let
$I=I(\mathcal{C})$ be 
its edge ideal. We prove that if ${\cal C}$
satisfies the packing property $($resp. max-flow min-cut property$)$,
then there is a  
uniform Cohen-Macaulay clutter ${\cal C}_1$ satisfying the packing
property $($resp. max-flow min-cut property$)$ such that ${\cal C}$
is a minor of  
${\cal C}_1$. For arbitrary edge ideals of clutters we prove that 
the normality property is closed under parallelizations. Then we show
some applications to edge ideals and clutters which are related to a 
conjecture of Conforti and Cornu\'ejols and to max-flow min-cut
problems. 
\end{abstract}

\section{Introduction}

Let $R=K[x_1,\ldots,x_n]$ be a polynomial ring 
over a field $K$ and let $I$ be an ideal 
of $R$  minimally generated by a finite set 
$F=\{x^{v_1},\ldots,x^{v_q}\}$ of square-free monomials. As usual we
use the 
notation $x^a:=x_1^{a_1} \cdots x_n^{a_n}$, 
where $a=(a_1,\ldots,a_n)$ is in $\mathbb{N}^n$. The {\it support\/}
of a monomial $x^a$ is given by ${\rm supp}(x^a)= \{x_i\, |\,
a_i>0\}$. For technical reasons we 
shall assume that each variable $x_i$ occurs in at least one
monomial of $F$. 

A {\it clutter\/} with finite vertex set 
$X$ is a family of subsets of $X$, called edges, none 
of which is included in another. The set of vertices of a clutter
$\mathcal{C}$ is denoted by $V(\mathcal{C})$ and the set of edges of
$\mathcal{C}$ is denoted by $E(\mathcal{C})$. A clutter is called 
$d$-{\it uniform\/} if all its edges have exactly 
$d$ vertices. We associate to the 
ideal $I$ a {\it clutter\/} $\cal C$ by taking the set 
of indeterminates $X=\{x_1,\ldots,x_n\}$ as vertex set and 
$E=\{S_1,\ldots,S_q\}$ as edge set, where $S_k$ is the support of
$x^{v_k}$. The vector $v_k$ is called the {\it characteristic
vector\/} of $S_k$. The assignment $I\mapsto \cal C$ gives a natural
one to one 
correspondence between the family of square-free monomial ideals and
the family of clutters. The ideal $I$ is called the {\it edge
ideal\/} 
of $\cal C$. To stress the relationship 
between $I$ and $\cal C$ we will use the notation $I=I({\cal C})$. 
Edge ideals of graphs were introduced and studied in \cite{ITG,Vi2}. 
Edge ideals of clutters also correspond to simplicial complexes via
the 
Stanley-Reisner correspondence \cite{Stanley} and to facet ideals 
\cite{faridi,zheng}. The Cohen-Macaulay property of edge ideals has
been recently studied in
\cite{carra-ferrarello,faridi1,linearquotients,MRV,bipartite_scm}  
using a combinatorial approach based 
on the notions of shellability, linear quotients, unmixedness, 
acyclicity and transitivity of digraphs, and the K\"onig property.

The aim of this note is to study the behavior, under certain
operations, of various algebraic
and combinatorial optimization properties of edge
ideals and clutters such as the Cohen-Macaulay property, the
normality, the torsion freeness, the packing and the max-flow min-cut 
properties. The study of edge ideals from the combinatorial
optimization point of view was initiated in \cite{unimod,shiftcon}
and  
continued in
\cite{multical,normali,alexdual,reesclu,clutters,perfect}, 
see also \cite{HHTZ}. The
Cohen-Macaulay and normality properties are two of the most
interesting 
properties an 
edge ideal can have, see  \cite{BHer,faridi1,Stanley,monalg} and 
\cite{huneke-swanson-book,bookthree} respectively.

Recall that the {\it integral closure\/} of $I^i$, denoted by
$\overline{I^i}$, is the ideal of $R$ given by  
$$
\overline{I^i}=(\{x^a\in R\vert\, \exists\, p\geq 1;(x^a)^{p}\in
I^{pi}\}).
$$ 
An ideal $I$ is called {\it normal\/} if $I^i=\overline{I^i}$ for
all $i$. A clutter obtained from $\mathcal{C}$ by a sequence of
deletions and duplications of vertices is called a {\it
parallelization\/} of 
$\mathcal{C}$ and a
clutter obtained from $\mathcal{C}$ by a sequence of deletions and
contractions of vertices is called a {\it minor\/} of $\mathcal{C}$, see 
Section~\ref{duplications-section}. It is known that
the normality of $I(\mathcal{C})$ is closed under
minors \cite{normali}. One of our main results shows that the
normality of $I(\mathcal{C})$
is closed under parallelizations:

\medskip

\noindent {\bf Theorem~\ref{parall-normal}}{\it\ Let $\mathcal{C}$ be
a clutter 
and let $\mathcal{C}'$ be a parallelization of $\mathcal{C}$. 
If $I(\mathcal{C})$ is normal, then $I(\mathcal{C}')$ is normal. 
}

\medskip

The ideal $I=I(\mathcal{C})$ is called {\it normally
torsion free\/} if $I^i=I^{(i)}$ for all
$i$, where $I^{(i)}$ is the $i{\it th}$ symbolic power 
of $I$.
As an application we prove that if $I(\mathcal{C})$ is normally
torsion free and $\mathcal{C}'$ is a parallelization of $\mathcal{C}$,
then $I(\mathcal{C}')$ is normally torsion free
(Corollary~\ref{jan26-08}). Let $A$ be the incidence matrix of
$\mathcal{C}$, i.e., $A$ is the matrix with column vectors
$v_1,\ldots,v_q$. A clutter $\cal C$ satisfies the {\it max-flow
min-cut\/} 
(MFMC) property if both sides 
of the LP-duality equation
$$
{\rm min}\{\langle w,x\rangle \vert\, x\geq 0; xA\geq{\mathbf 1}\}=
{\rm max}\{\langle y,{\mathbf 1}\rangle \vert\, y\geq 0; Ay\leq w\} 
$$
have integral optimum solutions $x$ and $y$ for each non-negative 
integral vector $w$. A remarkable result of \cite{clutters} (cf.
\cite[Theorem~4.6]{reesclu}) shows that $I(\mathcal{C})$ is normally
torsion free 
if and only if $\mathcal{C}$ has the max-flow min-cut property. This
fact makes a strong 
connection between 
commutative algebra and combinatorial optimization. 
It is known \cite[Chapter~79]{Schr2} that
a clutter $\mathcal{C}$ satisfies the max-flow min-cut property if
and only if all parallelizations of the clutter $\mathcal{C}$ satisfy
the K\"onig 
property (see Definition~\ref{konigpropdef}). As another application
we give a proof of this fact using that the integrality of the
polyhedron $\{x\vert\, x\geq 0;xA\geq\mathbf{1}\}$ is closed under
parallelizations and minors and using that the normality of
$I(\mathcal{C})$ is preserved under parallelizations and minors
(Corollary~\ref{mfmc-iff-cw}).

A clutter $\mathcal{C}$ satisfies the {\it packing property\/} (PP for
short) if all minors of $\mathcal{C}$ satisfy the K\"onig property. We
say that a clutter $\mathcal{C}$ is {\it Cohen-Macaulay\/} if
$R/I(\mathcal{C})$ is a Cohen-Macaulay ring, see \cite{Mats}. 
The other main result of this note is:

\medskip

\noindent {\bf Theorem~\ref{pp-mfmc-cm}}{\it\ Let $\mathcal{C}$ be a
$d$-uniform clutter on the vertex set $X$. Let 
$$Y=\{y_{ij}\vert\, 1\leq i\leq n;\, 1\leq
j\leq d-1\}$$
 be a set of new variables, and let ${\cal C}'$ be 
the clutter with vertex set $V({\cal C}')=X\cup Y$ and 
edge set 
$$E({\cal C}')=E({\cal C})\cup\{\{x_1,y_{11},\ldots,y_{1(d-1)}\},\ldots,
\{x_n,y_{n1},\ldots,y_{n(d-1)}\}\}.$$
Then the edge ideal $I({\cal C}')$ is Cohen-Macaulay. 
If $\cal C$ satisfies  {\rm PP} {\rm(}resp. max-flow min-cut{\rm)}, 
then ${\cal C}'$ satisfies  {\rm PP} {\rm(}resp. max-flow 
min-cut{\rm)}.
}

\medskip

It is well known that if $\mathcal{C}$ satisfies the max-flow min-cut
property, then $\mathcal{C}$ satisfies the packing property
\cite{cornu-book} (see Corollary~\ref{aug25-06}). Conforti and
Cornu\'ejols \cite{CC} conjecture that the converse is also
true. Theorem~\ref{pp-mfmc-cm} is interesting because it says that for
uniform clutters it suffices to prove the conjecture for Cohen-Macaulay
clutters, which have a rich structure.  
The Conforti-Cornu\'ejols conjecture has been studied 
in \cite{mfmc,reesclu,clutters} using an algebraic approach based on 
certain algebraic properties of blowup algebras. 

\section{Normality is preserved under
parallelizations}\label{duplications-section} 

Let $\mathcal{C}$ be a clutter on the vertex set
$X=\{x_1,\ldots,x_n\}$ and let
$I=I(\mathcal{C})=(x^{v_1},\ldots,x^{v_q})$ be its edge
ideal. The {\it incidence matrix\/} of $\mathcal{C}$, denoted by
$A=(a_{ij})$, is
the $n\times q$ matrix 
whose $(i,j)$ entry is given by $a_{ij}=1$ if $x_i\in g_j$ and
$a_{ij}=0$ otherwise, where $g_1,\ldots,g_q$ are the edges of
$\mathcal{C}$. Notice that the column vectors of $A$ are
$v_1,\ldots,v_q$. 
Recall that the {\it Rees
algebra\/} of $I$ is given by: 
$$
R[It]:=R\oplus It\oplus\cdots\oplus I^{i}t^i\oplus\cdots
\subset R[t],
$$
where $t$ is a new variable. The Rees algebra of $I$ can be written as
$$
R[It]=K[\{x^at^b\vert\, (a,b)\in\mathbb{N}{\cal A}'\}]
$$
where ${\cal A}'=\{(v_1,1),\ldots,(v_q,1),e_1,\ldots,e_n\}$ and 
$\mathbb{N}{\cal A}'$ is the subsemigroup of 
$\mathbb{N}^{n+1}$ spanned by ${\cal A}'$. In other words $R[It]$ is
equal to $K[\mathbb{N}{\cal A}']$, the semigroup ring of
$\mathbb{N}{\cal A}'$, see \cite{gilmer}. On the other hand according to
\cite[Theorem~7.2.28]{monalg} the 
integral closure of $R[It]$ 
in its field of fractions can be expressed as
\begin{eqnarray}
\overline{R[It]}&=&K[\{x^at^b\vert\, (a,b)\in
\mathbb{Z}{\cal
A}'\cap \mathbb{R}_+{\cal A}'\}]\nonumber\\ 
&=&R\oplus
\overline{I}t\oplus\overline{I^2}t^2\oplus\cdots\oplus
 \overline{I^i}t^i\oplus\cdots,\nonumber\label{jun05-1-03}
\end{eqnarray}
where $\overline{I^i}$ is the integral closure of $I^i$,
$\mathbb{R}_+{\cal A}'$ is the cone  
spanned by ${\cal A}'$, and $\mathbb{Z}{\cal A}'$ is the subgroup
spanned by ${\cal A}'$.  Notice that $\mathbb{Z}{\cal
A}'=\mathbb{Z}^{n+1}$. Hence $R[It]$ is normal if and
only if  
any of the following two equivalent conditions hold: 
\begin{description}
\item{\rm (a)} $\mathbb{N}{\cal A}'=
\mathbb{Z}^{n+1}\cap\mathbb{R}_+{\cal A}'$. \vspace{-1mm}
\item{\rm (b)} $I^{i}=\overline{I^i}$ for all $i\geq 1$.
\end{description}
If the second condition holds we say that $I$ is a {\it normal} 
ideal.  

Let $\mathcal{C}$ be a clutter on the vertex set
$X=\{x_1,\ldots,x_n\}$ and let $x_i\in X$. Then {\it duplicating\/}
$x_i$ means extending $X$ by a new vertex $x_i'$ and replacing
$E(\mathcal{C})$ by
$$    
E(\mathcal{C})\cup\{(e\setminus\{x_i\})\cup\{x_i'\}\vert\, x_i\in e\in
E(\mathcal{C})\}.
$$
The {\it deletion\/} of $x_i$, denoted by
$\mathcal{C}\setminus\{x_i\}$, is the clutter formed from
$\mathcal{C}$ by deleting the vertex $x_i$ and all edges containing
$x_i$. A clutter obtained from $\mathcal{C}$ by a sequence of
deletions and 
duplications of vertices is called a {\it parallelization\/}. If 
$w=(w_i)$ is a vector in $\mathbb{N}^n$, we denote by $\mathcal{C}^w$
the clutter obtained from
$\mathcal{C}$ by deleting any vertex $x_i$ with $w_i=0$ and
duplicating $w_i-1$ times any vertex $x_i$ if $w_i\geq 1$. The map
$w\mapsto \mathcal{C}^w$ gives a one to one correspondence between
$\mathbb{N}^n$ and the parallelizations of $\mathcal{C}$. 

\begin{Example}\rm Let $G$ be the graph whose only edge is 
$\{x_1,x_2\}$ and let $w=(3,3)$. Then
$G^w=\mathcal{K}_{3,3}$ is the complete bipartite graph 
with bipartition $V_1=\{x_1,x_1^2,x_1^3\}$ and
$V_2=\{x_2,x_2^2,x_2^3\}$. Notice that $x_i^k$ is a vertex, i.e., $k$
is an index not an exponent.
\end{Example}

The following notion of minor comes from combinatorial 
optimization \cite{CC,Schr2} and it is not apparently related 
to the minors (subdeterminants) of $A$.  

\begin{Definition}\rm Let 
$X'=\{x_{i_1},\ldots,x_{i_r},x_{j_1},\ldots,x_{j_s}\}$ be a subset of $X$. 
A {\it minor\/} of $I$ is a proper ideal $I'$ of $R'=K[X\setminus
X']$ obtained from  $I$ by making 
$x_{i_k}=0$ and $x_{j_\ell}=1$ for all $k,\ell$. The ideal $I$ is
considered itself a minor.  
A {\it minor\/} of $\cal C$ 
is a clutter ${\cal C}'$ whose edge ideal is $I'$.  
\end{Definition}

Notice that the generators of $I'$ are obtained from the generators of
$I$ by making $x_{i_k}=0$ and $x_{j_\ell}=1$ for all $k,\ell$. This
means that $\mathcal{C}'$ is obtained from $\mathcal{C}$ by shrinking
some edges and deleting some other edges. Also notice 
that ${\cal C}'$ is obtained from $I'$ by considering the unique set 
of square-free monomials of $R'$ that minimally generate $I'$. If
$I'$ is the ideal obtained from $I$ by making $x_i=0$, then 
$I'=I(\mathcal{C}\setminus\{x_i\})$, i.e., making a variable equal to
zero corresponds to a deletion. If $I'$ is the minor obtained from $I$ by 
making $x_i=0$ for $1\leq i\leq r$ and $x_i=1$ for $r+1\leq i\leq s$,
then in algebraic terms $I'$ 
can be expressed as 
$$
(I\cap
K[x_{r+1},\ldots,x_n])_\mathfrak{p}=I'K[x_{r+1},\ldots,x_n]_\mathfrak{p}, 
$$
where $(I\cap
K[x_{r+1},\ldots,x_n])_\mathfrak{p}$ and
$K[x_{r+1},\ldots,x_n]_\mathfrak{p}$ are localizations at 
the prime ideal $\mathfrak{p}$ 
generated by the variables $x_{s+1},\ldots,x_n$. 

\medskip

It is known that the normality of $I(\mathcal{C})$ is closed under
minors \cite{normali}. A main
result of this section shows that the normality of $I(\mathcal{C})$
is closed under parallelizations.

\begin{Theorem}\label{parall-normal} Let $\mathcal{C}$ be a clutter
and let $\mathcal{C}'$ be a parallelization of $\mathcal{C}$. 
If $I(\mathcal{C})$ is normal, then $I(\mathcal{C}')$ is normal. 
\end{Theorem}

\demo From \cite{normali} we obtain that if $I(\mathcal{C})$ is normal
and $\mathcal{C}'$ is a minor of $\mathcal{C}$, then
$I(\mathcal{C}')$ is also 
normal. Thus we need only show that the normality of $I(\mathcal{C})$
is preserved when we duplicate a vertex of $\mathcal{C}$. 
Let $V(\mathcal{C})=\{x_2,\ldots,x_n\}$ be the vertex set of
$\mathcal{C}$ and let $\mathcal{C}'$ be the clutter obtained from
$\mathcal{C}$ by duplication of the vertex
$x_2$. We denote the duplication of $x_2$ by $x_1$.  
We may assume that
$$
I=I(\mathcal{C})=(x_2x^{w_1},\ldots,x_2x^{w_r},x^{w_{r+1}},
\ldots,x^{w_q}),$$ 
where $x^{w_i}\in K[x_3,\ldots,x_n]$ 
for all $i$. We must show that the ideal 
$$
I(\mathcal{C}')=I+(x_1x^{w_1},\ldots,x_1x^{w_r})
$$
is normal. Consider the sets
\begin{eqnarray*}
\mathcal{A}&=&\{e_2,\ldots,e_n,(0,1,w_1,1),\ldots,
(0,1,w_r,1),(0,0,w_{r+1},1),\ldots,(0,0,w_q,1)\},\\ 
\mathcal{A}'&=&\mathcal{A}\cup\{e_1,(1,0,w_1,1),\ldots,
(1,0,w_r,1)\}.
\end{eqnarray*}
By hypothesis $\mathbb{Z}^{n+1}\cap\mathbb{R}_+{\cal
A}=\mathbb{N}{\cal A}$. We must prove that 
$\mathbb{Z}^{n+1}\cap\mathbb{R}_+{\cal
A}'=\mathbb{N}{\cal A}'$. It suffices to show that the left hand side
is contained in the right hand side because the other inclusion always
holds. Take an integral vector $(a,b,c,d)$ in 
$\mathbb{R}_+{\cal A}'$, where $a,b,d\in\mathbb{Z}$ and
$c\in\mathbb{Z}^{n-2}$. Then
$$
(a,b,c,d)=\sum_{i=1}^r\alpha_i(0,1,w_i,1)+\sum_{i=r+1}^q\alpha_i(0,0,w_i,1)+
\sum_{i=1}^r\beta_i(1,0,w_i,1)+\sum_{i=1}^n\gamma_ie_i
$$
for some $\alpha_i$, $\beta_i$, $\gamma_i$ in $\mathbb{R}_+$.
Comparing entries one has
\begin{eqnarray*}
a&=&\beta_1+\cdots+\beta_r+\gamma_1,\\
b&=&\alpha_1+\cdots+\alpha_r+\gamma_2,\\
c&=&\sum_{i=1}^r(\alpha_i+\beta_i)w_i+\sum_{i=r+1}^q\alpha_iw_i
+\sum_{i=3}^n\gamma_ie_i,\\
d&=&\sum_{i=1}^r(\alpha_i+\beta_i)+\sum_{i=r+1}^q\alpha_i.
\end{eqnarray*}
Consequently we obtain the equality 
$$
(0,a+b,c,d)=\sum_{i=1}^r(\alpha_i+\beta_i)(0,1,w_i,1)+
\sum_{i=r+1}^q\alpha_i(0,0,w_i,1)
+(\gamma_1+\gamma_2)e_2+\sum_{i=3}^n\gamma_ie_i,
$$
that is, the vector $(0,a+b,c,d)$ is in $\mathbb{Z}^{n+1}\cap
\mathbb{R}_+\mathcal{A}=\mathbb{N}\mathcal{A}$. Thus there are
$\lambda_i$, $\mu_i$ in $\mathbb{N}$ such that. 
$$
(0,a+b,c,d)=\sum_{i=1}^r\mu_i(0,1,w_i,1)+\sum_{i=r+1}^q\mu_i(0,0,w_i,1)
+\sum_{i=2}^n\lambda_ie_i.
$$
Comparing entries we obtain the equalities
\begin{eqnarray*}
a+b&=&\mu_1+\cdots+\mu_r+\lambda_2,\\
c&=&\mu_1w_1+\cdots+\mu_qw_q+ \lambda_3e_3+\cdots+\lambda_ne_n,
\\
d&=&\mu_1+\cdots+\mu_q.
\end{eqnarray*}

Case (I): $b\leq\sum_{i=1}^r\mu_i$. If $b<\mu_1$, we set $b=\mu_1'$, 
$\mu_1'<\mu_1$, and define $\mu_1''=\mu_1-\mu_1'$. Otherwise pick
$s\geq 2$ such  
that 
$$
\mu_1+\cdots+\mu_{s-1}\leq b\leq \mu_1+\cdots+\mu_{s}
$$
Then $b=\mu_1+\cdots+\mu_{s-1}+\mu_s'$, where $ \mu_s'\leq \mu_s$.
Set $\mu_s''=\mu_s-\mu_s'$. Notice that 
\begin{eqnarray*}
a+b&=&\mu_1+\cdots+\mu_r+\lambda_2=a+\mu_1+\cdots+\mu_{s-1}+\mu_s',\\ 
a&=&\mu_s+\cdots+\mu_r+\lambda_2-\mu_s'=\mu_{s+1}+\cdots
+\mu_r+\mu_s''+\lambda_2.
\end{eqnarray*}
Then
\begin{eqnarray*}
(a,b,c,d)&=&\sum_{i=1}^{s-1}\mu_i(0,1,w_i,1)+\mu_s'(0,1,w_s,1)
+\sum_{i=r+1}^q\mu_i(0,0,w_i,1)\\
& &+\mu_s''(1,0,w_s,1)+\sum_{i=s+1}^r\mu_i(1,0,w_i,1)
+\lambda_2e_1+\sum_{i=3}^n\lambda_ie_i,
\end{eqnarray*}
that is, $(a,b,c,d)\in\mathbb{N}{\mathcal A}'$.
 
Case (II): $b>\sum_{i=1}^r\mu_i$. Then
$b=\sum_{i=1}^r\mu_i+\lambda_2'$. Since 
$$
a+b=\mu_1+\cdots+\mu_r+\lambda_2=a+\mu_1+\cdots+\mu_r
+\lambda_2'
$$
we get $a=\lambda_2-\lambda_2'$. In particular $\lambda_2\geq
\lambda_2'$. Then
$$
(a,b,c,d)=\sum_{i=1}^r\mu_i(0,1,w_i,1)+\sum_{i=r+1}^q\mu_i(0,0,w_i,1)
+ae_1+\lambda_2'e_2+\sum_{i=3}^n\lambda_ie_i
$$
that is, $(a,b,c,d)\in\mathbb{N}{\mathcal A}'$. \QED

\medskip

Our next goal is to present some applications of this result, but
first we need to prove a couple of lemmas and we need to 
recall some notions and results.

\begin{Definition}\rm A subset $C\subset X$ is a 
{\it minimal vertex cover\/} of the clutter $\cal C$ if: 
(i) every edge of $\cal C$ contains at least one vertex of $C$, 
and (ii) there is no proper subset of $C$ with the first 
property. If $C$ satisfies condition (i) only, then $C$ is 
called a {\it vertex cover\/} of $\cal C$.   
\end{Definition}

\begin{Definition}\rm Let $A$ be the incidence matrix of
$\mathcal{C}$. The clutter $\cal C$ satisfies the {\it max-flow min-cut\/}
(MFMC) 
property if both sides 
of the LP-duality equation
\begin{equation}\label{jun6-2-03-1}
{\rm min}\{\langle w,x\rangle \vert\, x\geq 0; xA\geq{\mathbf 1}\}=
{\rm max}\{\langle y,{\mathbf 1}\rangle \vert\, y\geq 0; Ay\leq w\} 
\end{equation}
have integral optimum solutions $x$ and $y$ for each non-negative 
integral vector $w$. 
\end{Definition}

Let $A$ be the {\it incidence matrix\/} of $\cal C$ whose column 
vectors are $v_1,\ldots,v_q$. The {\it set covering polyhedron\/} 
 of ${\cal C}$ is given by:
$$
Q(A)=\{x\in\mathbb{R}^n\vert\, x\geq 0;\, xA\geq{\mathbf 1}\},
$$
where $\mathbf{1}=(1,\ldots,1)$. This polyhedron was studied in
\cite{reesclu,clutters} to  
characterize the max-flow min-cut property of $\cal C$ and to study 
certain algebraic properties of blowup algebras. A clutter $\cal C$
is said  
to be {\em ideal} if $Q(A)$ is an {\it integral polyhedron\/}, i.e.,
it has 
only integral vertices.  
The integral vertices of $Q(A)$ are precisely the characteristic
vectors 
of the minimal vertex covers of $\cal C$
\cite[Proposition~2.2]{reesclu}.

\begin{Theorem}[\rm\cite{normali,reesclu,clutters,HuSV}]\label{noclu1}
The following conditions are equivalent{\rm :} 
\begin{description}
\item{\rm(i)\ \,} ${\rm gr}_I(R)=R[It]/IR[It]$ is reduced, i.e.,
${\rm gr}_I(R)$ has no non-zero nilpotent elements.
\vspace{-1mm}
\item{\rm (ii)\ } $R[It]$ is normal and $Q(A)$ is an integral
polyhedron.\vspace{-1mm}
\item{\rm (iii)\,}  $I^{i}=I^{(i)}$ for $i\geq 1$, where $I^{(i)}$ is
the $i${\it th} symbolic power of $I$.
\vspace{-1mm}
\item{\rm (iv)\ } $\cal C$ has the max-flow min-cut property.
\end{description}
\end{Theorem}

If condition (iii) is satisfied we say that $I$ is {\it normally torsion
free\/}. A set of edges of the clutter $\cal C$ is {\it
independent\/} or {\it stable} if no two of them have a
common vertex. 
We denote the smallest number of vertices in any 
minimal vertex cover of $\cal C$ by $\alpha_0({\cal C})$ and the 
maximum number of independent edges of ${\cal C}$ by $\beta_1({\cal
C})$. These  
numbers are related to min-max problems because they satisfy:
\begin{eqnarray*}
\lefteqn{\alpha_0({\cal C})\geq {\rm min}\{\langle {\mathbf 1},x\rangle \vert\, 
x\geq 0; xA\geq {\mathbf 1}\}}\\
&\ \ \ \ \ \ \ \ \ \ &
={\rm max}\{\langle y,{\mathbf 1}\rangle \vert\, y\geq 0; Ay\leq{\mathbf 1}\}
\geq \beta_1({\cal C}). 
\end{eqnarray*}
Notice that $\alpha_0({\cal C})=\beta_1({\cal C})$ if and only if
both sides of  
the equality have integral optimum solutions. These two numbers can
be interpreted 
in terms of invariants of $I$. By \cite{reesclu} the height 
of the ideal $I$, denoted by ${\rm ht}(I)$, is equal to the {\it
vertex covering number\/} $\alpha_0({\cal C})$ and the 
{\it edge independence number\/} $\beta_1({\cal C})$ is equal to the 
maximum $r$ such that there exists a regular sequence of $r$ monomials
inside $I$.  

\begin{Definition}\label{konigpropdef}\rm If $\alpha_0({\cal
C})=\beta_1({\cal C})$ we say  
that the clutter $\cal C$ (or the ideal $I$) has the {\it K\"onig
property\/}.
\end{Definition}

\begin{Definition}\rm The clutter $\cal C$ (or the ideal $I$) satisfy
the {\it packing property\/}  
(PP for short) if all its minors satisfy the K\"onig property, i.e., 
$\alpha_0({\cal C}')=\beta_1({\cal C}')$ 
for any minor ${\cal C}'$ of $\cal C$.
\end{Definition}

\begin{Theorem}{\rm(A. Lehman; see \cite[Theorem~1.8]{cornu-book})}
\label{lehman} If\, $\cal C$ has the packing property, then $Q(A)$ is
integral.
\end{Theorem}

\begin{Corollary}[\rm\cite{cornu-book}]\label{aug25-06} If the clutter
$\cal C$ has the max-flow
min-cut property,  
then $\cal C$ has the packing property. 
\end{Corollary}

\demo Assume that the clutter $\mathcal{C}$ has the max-flow min-cut
property. This property is closed under taking minors. Thus it 
suffices to prove that $\mathcal C$ has 
the K\"onig property. We denote the incidence matrix of ${\cal C}$ by
$A$. By hypothesis the LP-duality equation
$$
{\rm min}\{\langle\mathbf{1},x\rangle \vert\, x\geq 0; xA\geq \mathbf{1}\}=
{\rm max}\{\langle y,\mathbf{1}\rangle \vert\, y\geq 0; Ay\leq\mathbf{1}\}
$$
has optimum integral solutions $x$, $y$. To complete the 
proof notice that the left hand side of this 
equality is $\alpha_0({\cal C})$ and the right hand side 
is $\beta_1({\cal C})$. \QED

\medskip

Conforti and Cornu\'ejols conjecture that the converse is also
true:

\begin{Conjecture}{\rm(\cite{CC})}\label{conforti-cornuejols1}
\rm  If the clutter $\cal C$ has the
packing property, then $\cal C$ has the max-flow min-cut property.
\end{Conjecture} 

To the best of our knowledge this conjecture is open, 
see \cite[Conjecture~1.6]{cornu-book}.

\begin{Corollary}\label{jan26-08} Let $\mathcal{C}$ be a clutter and 
let $\mathcal{C}'$ be a parallelization of $\mathcal{C}$. 
If $I(\mathcal{C})$ is normally torsion free, 
then $I(\mathcal{C}')$ is normally torsion free.
\end{Corollary}

\demo Let $A$ and $A'$ be the incidence matrices of $\mathcal{C}$ and
$\mathcal{C}'$ respectively. 
By Theorem~\ref{noclu1} the ideal $I(\mathcal{C})$ is normal and
$Q(A)$ is integral. From Theorem~\ref{parall-normal} the ideal
$I(\mathcal{C}')$ is normal, and since the integrality of $Q(A)$ is 
closed under minors and parallelizations (see \cite{reesclu} and 
\cite{Schr2}) we get that $Q(A')$ is again integral. Thus applying 
Theorem~\ref{noclu1} once more we get that $I(\mathcal{C}')$ is
normally torsion free. \QED

\begin{Corollary}\label{mfmc-parall-konig} Let $\mathcal{C}$ be a
clutter and  
let $\mathcal{C}'$ be a parallelization of $\mathcal{C}$. 
If $\mathcal{C}$ has the max-flow min-cut property, then
$\mathcal{C}'$ has the K\"onig property. In particular 
$\mathcal{C}^w$ has the K\"onig property for all $w\in\mathbb{N}^n$. 

\end{Corollary}

\demo By Corollary~\ref{jan26-08} the clutter $\mathcal{C}'$ has the
max-flow min-cut property. Thus applying Corollary~\ref{aug25-06} we
obtain that $\mathcal{C}'$ has the K\"onig property. \QED

\begin{Lemma}\label{beta-w} Let $\mathcal{C}$ be a clutter and let
$A$ be its 
incidence matrix. If $w=(w_i)$ is a vector in $\mathbb{N}^n$, then
$$
\beta_1(\mathcal{C}^w)\leq \max\{\langle y,\mathbf{1}\rangle\vert\, y\in
\mathbb{N}^q;\, Ay\leq w\}.
$$
\end{Lemma}

\demo We may assume that $w=(w_1,\ldots,w_m,0,\ldots,0)$, where
$w_i\geq 1$ for $i=1,\ldots,m$. Recall that for each $i$ the vertex $x_i$ is 
duplicated $w_i-1$ times. We denote the duplications of $x_i$ by 
$x_{i}^2,\ldots,x_i^{w_i}$ and set $x_i^1=x_i$. Thus the vertex set of 
$\mathcal{C}^w$ is equal to
$$
V(\mathcal{C}^w)=\{x_1^1,\ldots,x_1^{w_1},\ldots,x_i^1,\ldots,x_i^{w_i},
\ldots,x_m^1,\ldots,x_m^{w_m}\}.
$$
There are $f_1,\ldots,f_{\beta_1}$ independent edges of $\mathcal{C}^w$, where 
$\beta_1=\beta_1(\mathcal{C}^w)$. Each $f_i$ has the form
$$
f_k=\{x_{k_1}^{j_{k_1}},x_{k_2}^{j_{k_2}},\ldots,x_{k_r}^{j_{k_r}}\}\ \
\ \ \ 
(1\leq k_1<\cdots<k_r\leq m;\ 1\leq j_{k_i}\leq w_{k_i}).
$$
We set $g_k=\{x_{k_1}^{1},x_{k_2}^{1},\ldots,x_{k_r}^{1}\}=
\{x_{k_1},x_{k_2},\ldots,x_{k_r}\}$. By definition of $\mathcal{C}^w$
we get that $g_k\in E(\mathcal{C})$ for all $k$. We may re-order the $f_i$ so that 
$$
\underbrace{g_1=g_2=\cdots=g_{s_1}}_{s_1},\underbrace{g_{s_1+1}=\cdots=
g_{s_2}}_{s_2-s_1},\ldots,\underbrace{g_{s_{r-1}+1}=\cdots=g_{s_r}}_{s_r-s_{r-1}}
$$
and $g_{s_1},\ldots,g_{s_r}$ distinct, where $s_r=\beta_1$. Let
$v_{i}$ be the characteristic vector of $g_{s_i}$. 
Set $y=s_1e_{1}+(s_2-s_1)e_{2}+\cdots+(s_r-s_{r-1})e_{r}$. We may
assume that the incidence matrix $A$ of $\mathcal{C}$ has column
vector $v_1,\ldots,v_q$.   
Then $y$ satisfies $\langle y,\mathbf{1}\rangle=\beta_1$. For each $k_i$ the
number of variables of the form $x_{k_i}^\ell$ that occur in
$f_1,\ldots,f_{\beta_1}$ is at most $w_{k_i}$ because the $f_i$ are
pairwise disjoint. Hence for each $k_i$ the number of times that the
variable $x_{k_i}^1$ occurs in $g_1,\ldots,g_{\beta_1}$ is at most
$w_{k_i}$. Then 
$$
Ay=s_1v_{1}+(s_2-s_1)v_{2}+\cdots+(s_r-s_{r-1})v_{r}\leq w.
$$
Therefore we obtain the required inequality. \QED

\medskip

\medskip Let $\mathcal{C}$ be a clutter. For use below we denote the
set of minimal vertex covers of  
$\mathcal{C}$ by $\Upsilon(\mathcal{C})$.

\begin{Lemma}\label{alpha-w} Let $\mathcal{C}$ be a clutter and let $A$ be its
incidence matrix. If $w=(w_i)$ is a vector in $\mathbb{N}^n$, then
$$
\left.\min\left\{\sum_{x_i\in
C}w_i\right\vert\,
C\in\Upsilon(\mathcal{C})\right\}
= \alpha_0(\mathcal{C}^w).
$$
\end{Lemma}

\demo We may assume that
$w=(w_1,\ldots,w_m,w_{m+1},\ldots,w_{m_1},0,\ldots,0)$, where $w_i\geq
2$ for $i=1,\ldots,m$, $w_i=1$ 
for $i=m+1,\ldots,m_1$, and $w_i=0$ for $i>m_1$. 
Thus for $i=1,\ldots,m$ the vertex $x_i$ is
duplicated $w_i-1$ times. We denote the duplications of $x_i$ by 
$x_{i}^2,\ldots,x_i^{w_i}$ and set $x_i^1=x_i$. 

We first prove that the left hand side is less or equal than the right hand
side. Let $C$ be a minimal
vertex cover of $\mathcal{C}^w$ with $\alpha_0$ elements, 
where $\alpha_0=\alpha_0(\mathcal{C}^w)$. We may assume that 
$C\cap \{x_1,\ldots,x_{m_1}\}=\{x_1,\ldots,x_s\}$. Note that 
$x_i^1,\ldots,x_i^{w_i}$ are in $C$ for $i=1,\ldots,s$. Indeed 
since $C$ is a minimal vertex cover of $\mathcal{C}^w$, there exists
an edge $e$ of $\mathcal{C}^w$ such that $e\cap C=\{x_i^1\}$. Then
$(e\setminus\{x_i^1\})\cup\{x_i^j\}$ is an edge of $\mathcal{C}^w$ for
$j=1,\ldots, w_i$. Consequently  $x_i^j\in C$ for $j=1,\ldots,w_i$. Hence 
\begin{equation}\label{jan1-08}
w_1+\cdots+w_s\leq |C|=\alpha_0.
\end{equation}
On the other hand the set 
$C'=\{x_1,\ldots,x_s\}\cup\{x_{m_1+1},\ldots,x_n\}$ 
is a vertex cover of $\mathcal{C}$. Let $D$ be a minimal vertex cover
of $\mathcal{C}$ contained in $C'$. Let $e_D$ denote the characteristic
vector of $D$. Then, since $w_i=0$ for $i>m_1$, using
Eq.~(\ref{jan1-08}) we get
$$
\langle w,e_D\rangle=\sum_{x_i\in D}w_i=\sum_{x_i\in
D\cap\{x_1,\ldots,x_s\}}\hspace{-5mm}w_i\ \ \ \ \leq
\sum_{x_i\in\{x_1,\ldots,x_s\}}\hspace{-5mm}w_i\ \leq\ \alpha_0.
$$
This completes the proof of the asserted inequality. 

Next we show that the right hand side of the inequality 
is less or equal than the left hand side. Let $C$ be a minimal vertex
cover of $\mathcal{C}$. Note that the set
$$
C'=\cup_{x_i\in C}\{x_i^1,\ldots,x_i^{w_i}\}
$$
is a vertex cover of $\mathcal{C}^w$. Indeed any edge $e^w$ of
$\mathcal{C}^w$ has the form $e^w=\{x_{i_1}^{j_1},\ldots,x_{i_r}^{j_r}\}$
for some edge $e=\{x_{i_1},\ldots,x_{i_r}\}$ of $\mathcal{C}$ and
since $e$ is covered by $C$, we have that $e^w$ is covered by $C'$. 
Hence $\alpha_0(\mathcal{C}^w)\leq |C'|=\sum_{x_i\in C}w_i$. As $C$ was an
arbitrary vertex cover of $\mathcal{C}$ we get the asserted 
inequality. \QED

\begin{Corollary}{\rm\cite[Chapter~79]{Schr2}}\label{mfmc-iff-cw} Let $\mathcal{C}$ be
a clutter. Then $\mathcal{C}$ 
satisfies the max-flow min-cut property if and only if 
$\beta_1(\mathcal{C}^w)=\alpha_0(\mathcal{C}^w)$ for all
$w\in\mathbb{N}^n$.
\end{Corollary}

\demo If $\mathcal{C}$ has the max-flow min-cut property, 
then $\mathcal{C}^w$ has the K\"onig property by
Corollary~\ref{mfmc-parall-konig}. Conversely if $\mathcal{C}^w$ has
the K\"onig property for all $w\in\mathbb{N}^n$, then by 
Lemmas~\ref{beta-w} and \ref{alpha-w} both sides of the LP-duality
equation
$$
{\rm min}\{\langle w,x\rangle \vert\, x\geq 0; xA\geq{\mathbf 1}\}=
{\rm max}\{\langle y,{\mathbf 1}\rangle \vert\, y\geq 0; Ay\leq w\} 
$$
have integral optimum solutions $x$ and $y$ for each non-negative 
integral vector $w$, i.e., $\mathcal{C}$ has the max-flow min-cut
property. \QED

\section{Cohen-Macaulay ideals with max-flow min-cut}

One of the aims here is to show how to construct Cohen-Macaulay clutters
satisfying max-flow min-cut, PP, and normality properties. Let
$\mathcal{C}$ be a uniform clutter. A main result of this section
proves that if ${\cal C}$
satisfies {\rm PP} $($resp. max-flow min-cut$)$, then there is a 
uniform Cohen-Macaulay clutter ${\cal C}_1$ satisfying {\rm PP}
$($resp. max-flow min-cut$)$ such that ${\cal C}$ is a minor of 
${\cal C}_1$. In particular for uniform clutters we prove
that it suffices to show Conjecture~\ref{conforti-cornuejols1} for
Cohen-Macaulay clutters  (see Corollary~\ref{pp-cm}).
  
Let $R=K[x_1,\ldots,x_n]$ be a polynomial ring over a field $K$ and
let $\mathcal{C}$ be a clutter on the vertex set $X$. As usual, in
what follows, we
denote the edge ideal of $\mathcal{C}$ by $I=I(\mathcal{C})$. 
Recall that 
$\mathfrak{p}$ is a minimal prime of $I =I(\mathcal{C})$ if and only if 
$\mathfrak{p}=(C)$ for some minimal vertex cover $C$ of $\mathcal{C}$ 
\cite[Proposition~6.1.16]{monalg}. Thus the primary decomposition of
the edge ideal of $\mathcal{C}$ is given by
$$
I(\mathcal{C})=(C_1)\cap (C_2)\cap\cdots\cap (C_p),
$$
where $C_1,\ldots,C_p$ are the minimal vertex 
covers of $\mathcal{C}$. In particular observe that the {\it
height\/} 
of $I(\mathcal{C})$, denoted by ${\rm ht}\, I(\mathcal{C})$, is equal
to the minimum cardinality of a minimal vertex cover of 
$\mathcal{C}$. Also notice that the associated primes 
of $I(\mathcal{C})$ are precisely the minimal primes of
$I(\mathcal{C})$. 

\begin{Proposition}\label{may30-06}
Let $R[z_1,\ldots,z_\ell]$ be a polynomial ring over $R$. If $I$ is a 
normal ideal of $R$, then $J=(I,x_1z_1\cdots z_\ell)$ is a 
normal ideal of $R[z_1,\ldots,z_\ell]$.
\end{Proposition}
\demo By induction on $p$ we will show $\overline{J^p}=J^p$ for all $p\geq 1$. If
$p=1$, then $\overline{J}=J$ because $J$ is square-free (see
\cite[Corollary~7.3.15]{monalg}). Assume
$\overline{J^i}=J^i$ for  $i<p$ and $p\geq 2$. Let $y$ be a monomial in
$\overline{J^p}$, then $y^m\in J^{pm}$, for some $m>0$. Since $\overline{J^p}\subset 
\overline{J^{p-1}}=J^{p-1}$ we can write 
$$y=z_1^{t_1}\cdots z_\ell^{t_\ell}(x_1z_1\cdots z_\ell)^rMf_{1}\cdots
f_{p-r-1},$$
where $M$ is a monomial with $z_i\notin
{\rm supp}(M)$ for all $i$ and the $f_i$\/'s are 
monomials in $J$ with $z_i\notin{\rm supp}(f_j)$ for all
$i,j$. We set $h=Mf_{1}\cdots
f_{p-r-1}$. It suffices to 
show that $y\in J^p$. Since $y^m\in
J^{pm}$ we have  
\begin{equation}\label{powerm1}
y^m=z_1^{mt_1}\cdots z_\ell^{mt_\ell}(x_1z_1\cdots z_\ell)^{rm}h^m
=N(x_1z_1\cdots z_\ell)^sg_1\cdots g_{mp-s},
\end{equation} 
where $N$ is a monomial, $z_i\notin{\rm supp}(g_j)$ for all
$i,j$, and the $g_i$\/'s are monomials in $J$. We distinguish 
two cases: 

Case (a): Assume $t_i=0$ for some $i$, then $s\leq rm$ because $z_i^{rm}$ is
the maximum power of $z_i$ that divides $y^m$. Making $z_j=1$ for 
$j=1,\ldots,\ell$ in Eq.~(\ref{powerm1}) we get
$$
x_1^{rm-s}h^m
=N'g_1\cdots g_{mp-s}.
$$
Thus $h^m\in I^{(mp-s)-(rm-s)}=I^{m(p-r)}$. 
Therefore we get $h\in \overline{I^{p-r}}=I^{p-r}$ and 
$y=z_1^{t_1}\cdots z_\ell^{t_\ell}(x_1z_1\cdots z_\ell)^rh\in J^p$.  

Case (b): If $t_i>0$ for all $i$, we may assume $x_1\notin {\rm supp}(M)$, otherwise
$y\in J^p$. We may also assume $x_1\notin {\rm supp}(f_i)$ for all
$i$, otherwise it is not hard to see that we are back in case 
(a). Notice that $s\leq rm$, because $x_1\notin {\rm supp}(h)$. From 
Eq.~(\ref{powerm1}) it follows that $h\in 
\overline{I^{p-r}}=I^{p-r}$ and $y=z_1^{t_1}\cdots
z_\ell^{t_\ell}(x_1z_1\cdots z_\ell)^rh\in J^p$. \QED

\begin{Lemma}\label{may29-06} Let $R[z_1,\ldots,z_\ell]$ be a polynomial ring over $R$
and let $I_1$ be the ideal obtained from $I$ by 
making $x_1=0$. Then: {\rm (a)} if $I$ and $I_1$ satisfy the K\"onig property, then 
the ideal $J=(I,x_1z_1\cdots z_\ell)$ satisfies the K\"onig property,
and {\rm (b)} if $I$ satisfies {\rm PP}, then $J$ satisfies {\rm PP}.
\end{Lemma}

\demo (a): If ${\rm ht}(I)={\rm ht}(J)$, then $J$ satisfies
K\"onig because $I$ does. Assume that $g={\rm ht}(I)<{\rm ht}(J)$. Then 
${\rm ht}(J)=g+1$. Notice that every associated prime ideal of $I$ of
height $g$ cannot contain $x_1$. We claim that ${\rm ht}(I_1)=g$. If
$r={\rm ht}(I_1)<g$, pick a minimal prime $\mathfrak{p}$ of $I_1$ of
height $r$. Then $\mathfrak{p}+(x_1)$ is a prime ideal of height at 
most $g$ containing both $I$ and $x_1$, a contradiction. This proves
the claim. Since $I_1$ satisfies K\"onig, there are $g$ independent 
monomials in $I_1$. Hence $h_1,\ldots,h_g,x_1z_1\cdots z_\ell$ are
$g+1$ independent monomials in $J$, as required. Part (b) 
follows readily from part (a). \QED

\begin{Theorem}\label{pp-mfmc-cm} 
Let $\mathcal{C}$ be a
$d$-uniform clutter on the vertex set $X$. Let 
$$Y=\{y_{ij}\vert\, 1\leq i\leq n;\, 1\leq
j\leq d-1\}$$
 be a set of new variables, and let ${\cal C}'$ be 
the clutter with vertex set $V({\cal C}')=X\cup Y$ and 
edge set 
$$E({\cal C}')=E({\cal C})\cup\{\{x_1,y_{11},\ldots,y_{1(d-1)}\},\ldots,
\{x_n,y_{n1},\ldots,y_{n(d-1)}\}\}.$$
Then the edge ideal $I({\cal C}')$ is Cohen-Macaulay. 
If $\cal C$ satisfies  {\rm PP} {\rm(}resp. max-flow min-cut{\rm)}, 
then ${\cal C}'$ satisfies  {\rm PP} {\rm(}resp. max-flow 
min-cut{\rm)}.
\end{Theorem}

\demo Set $S=K[X\cup Y]$ and $I'=I({\cal C}')$. The clutter
$\mathcal{C}'$ is a grafting of $\mathcal{C}$ as defined by Faridi in
\cite{faridi1}. Then $I'$ is Cohen-Macaulay by 
\cite[Theorem~8.2]{faridi1}. If $\cal C$ satisfies PP,
then from Lemma~\ref{may29-06}(b) it follows that ${\cal C}'$
satisfies PP. Assume that $\cal C$ satisfies MFMC. By
Proposition~\ref{may30-06} $S[I't]$ is normal. Since ${\cal C}'$ 
satisfies PP, by Lehman's theorem we get that $Q(A')$ is integral,
where $A'$ is the incidence matrix of ${\cal C}'$. Therefore 
using Theorem~\ref{noclu1} we conclude that ${\cal C}'$ has MFMC. \QED 

\medskip

Recall that a clutter $\cal C$ is called {\it Cohen-Macaulay\/} (CM
for short) if $R/I({\cal C})$ is a Cohen-Macaulay ring. Since
$\mathcal{C}$ is a minor of $\mathcal{C}'$ we obtain:

\begin{Corollary}\label{pp-cm} Let $\cal C$ be a uniform clutter. If ${\cal C}$
satisfies {\rm PP} $($resp. max-flow min-cut$)$, then there is a 
uniform Cohen-Macaulay clutter ${\cal C}_1$ satisfying {\rm PP}
$($resp. max-flow min-cut$)$ such that ${\cal C}$ is a minor of 
${\cal C}_1$. 
\end{Corollary}

This result is interesting because it says that for 
uniform clutters it suffices to prove
Conjecture~\ref{conforti-cornuejols1} 
for Cohen-Macaulay
clutters, which have a rich structure.

\bibliographystyle{plain}

\end{document}